\documentclass[12pt]{article}
\usepackage{amsfonts}
\usepackage{amsmath}
\usepackage{}

\usepackage{amssymb,amsthm,amsmath,hyperref,txfonts}

\usepackage{graphicx,float}

\numberwithin{equation}{section}

\newtheorem{thm}{Theorem}[section]

\newtheorem{rem}{Remark}[section]

\theoremstyle{definition}

\allowdisplaybreaks

\begin{document}
\title{\bf The  Limit Behavior of the Second-Grade Fluid Equations}
\author{
\begin{tabular}{cc}
&  Wenhuo Su$^{1}$\thanks{Corresponding
author. E-mail:
suwenhuo@jxycu.edu.cn}, \ \ \ \ Aibin Zang$^1$\thanks{E-mail:
zangab05@126.com}\\
& \textit{\small 1.  Center of Applied Mathematics, Yichun University, Yichun 336000, P.R. China}\\[1.5mm]
\end{tabular}
}


\maketitle
\begin{abstract}
In the paper,  the limit behavior of solutions to the second-grade fluid system with no-slip boundary conditions is studied as both $\nu \rightarrow 0$ and  $\alpha \rightarrow 0$.
More precisely, it is verified that the convergence from second-grade fluid system to Euler system holds as $\nu \rightarrow 0$ and  $\alpha \rightarrow 0$ independently under the radial symmetry case.

\vspace{4mm}

 {\textbf{Keyword:}} Second-Grade  fluid; Euler equation; Bessel equation; Vanishing Viscosity limit; Radially symmetry  \\

\end{abstract}

\section{Introduction}

The 2D second-grade fluid equations are a viscous, incompressible fluid system given by
\begin{equation}\label{1.01}
  \partial_t(u - \alpha^2u) - \nu\Delta u + (u \cdot \nabla)(u - \alpha^2u) +  \sum_{j = 1}^2(u - \alpha^2u)_j\nabla u_j + \nabla P = 0,
\end{equation}
with incompressibility condition
\begin{equation}\label{1.02}
  \nabla \cdot u = 0.
\end{equation}
Where $u = (u_1, u_2)$ is the velocity and $P$ is the scalar pressure. The second-grade fluid is a model of non-Newtonian fluids described in \cite{Cioranescu1984Existence}, as one among a hierarchy of models of viscoelastic fluids, with
two parameters: $\alpha > 0$ corresponding to the elastic response, and $\nu > 0$ presenting viscosity.

In 2D smooth simply connected domain $\Omega$, we impose the boundary conditions are that no-slip boundary conditions i.e.
\begin{equation}\label{1.04}
  u(x, t) = 0,   \ \ \ \ \ \ \ \ \ \ \text{on} \ \ \partial\Omega.
\end{equation}
The existence and uniqueness of  the system (\ref{1.01}) and (\ref{1.02}) with boundary conditions (\ref{1.04}) were obtained by some mathematicians and we  can refer  to \cite{Cioranescu1984Existence}, \cite{Marsden2000The} and \cite{Analysis}.
We firstly assume  $\alpha = 0$ in (\ref{1.01}), then formally one obtains the 2D incompressible Navier-Stokes equations
\begin{equation}\label{1.03}
  \partial_tu^{NS} - \nu\Delta u^{NS} + u^{NS}\cdot \nabla u^{NS} + \nabla P^{NS} = 0, \ \ \ \nabla\cdot u^{NS} = 0.
\end{equation}
The term $\sum_{j = 1}^2u_j\nabla u_j = \nabla\big(\frac{1}{2}|u|^2\big)$ in (\ref{1.01}) can be absorbed by the pressure. The limit behaviors as $\alpha \rightarrow 0$ of second-grade fluid (\ref{1.01}) to the Navier-Stokes equations (\ref{1.03}) were studied in \cite{Busuioc2018From}, \cite{Busuioc1999On} and \cite{YanpingCao2009On}.
Secondly,  it is well known that the inviscid second-grade fluid are called Euler-$\alpha$ equations. According to the uniform estimates with viscosity in \cite{Cioranescu1984Existence}, Busuioc verified that vanishing viscosity limit to the Euler$-\alpha$ equations in $H^s$ in \cite{Busuioc1999On}. The limit of the model (\ref{1.01}) as either $\alpha \rightarrow 0$ or $ \nu \rightarrow 0$, is not related to the boundary layer, due to the same boundary conditions.
It is interesting to study the independent limits of second-grade fluid (\ref{1.01}),
 as $\alpha \rightarrow 0, \nu \rightarrow 0$  with vary kinds of boundary  condition. In \cite{Busuioc2012Incompressible}, Busuioc et.al.  obtained the limit of second-grade fluid equations  under Navier boundary conditions, as $\alpha \rightarrow 0, \nu \rightarrow 0$, independently  is described by the incompressible Euler equations

\begin{equation}\label{1.05}
   \partial_t{u}^E + {u}^E\cdot\nabla{u}^E + \nabla{{P}^E} = 0, ~\nabla\cdot u^{E}=0,
\end{equation}
with no-penetration condition ${u}^E\cdot n = 0$.

The second author and his collaborators in \cite{Lopes2015Convergence} proved that the convergence from  Euler-$\alpha$ equations  to Euler equations (\ref{1.05})  with no-slip boundary conditions  as $\alpha\to 0$. For no-slip boundary conditions, the convergence of second-grade fluid to Euler equation is more complicated, as $\alpha \rightarrow 0,\nu \rightarrow 0$, simultaneously.  Since we have obtained the limit behavior from second-grade fluid to Navier-Stokes in \cite{Busuioc1999On} and \cite{Busuioc2018From}, however it is a challenging
problem for the vanishing viscosity limit of  Navier-Stokes solutions near a solid boundary.

In \cite{Filho2015Approximation}, the second author and his collaborators studied the limit behavior of the solution $u^{(\alpha, \nu)}$ to (\ref{1.01}),(\ref{1.02}) with (\ref{1.04}) as both $\alpha$ and $\nu$ go to zero.
The results above can be illustrated by the following figure.
Under suitable initial velocity assumptions, we observe that if the viscosity $\nu = \mathcal {O}(\alpha^2) $, as $\alpha \rightarrow 0$, as shown in region \uppercase\expandafter{\romannumeral4},  the solution $u^{(\alpha, \nu)}$ converges  to the solution $u^E$ of Euler equations in $L^2$ space uniformly in time.  When $\alpha^2 << \nu = \mathcal {O}(\alpha^{6/5})$ lying in  region \uppercase\expandafter{\romannumeral3} and  $\alpha = \mathcal {O}(\nu^{\frac{1}{2}}) $ filling in region \uppercase\expandafter{\romannumeral1}, one can obtains the necessary and sufficient conditions for the convergence are Kato's conditions with the regions near the boundary width of $\mathcal {O}(\alpha^3\nu^{-3/2})$ and $\mathcal {O}(\nu)$, respectively.  However, we have no results in region \uppercase\expandafter{\romannumeral2}.


\centerline{\includegraphics[width=6cm, height=6cm]{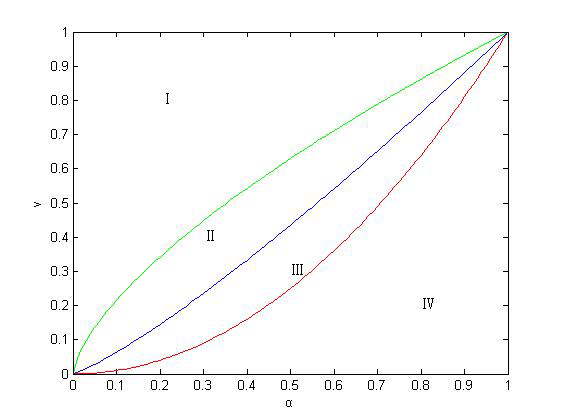}}

In this paper, we further discuss the limit behavior as $\nu, \alpha \rightarrow 0$ of the solution to the problem (\ref{1.01}),(\ref{1.02}) and (\ref{1.04}). Especially, we can check the convergence in region
\uppercase\expandafter{\romannumeral2} under the radial symmetry case. To verify this problem,  firstly we  take the curl operator to the equation (\ref{1.01}) and set $\omega = curl u$ the vorticity of $u$. Thus one has the 2D vorticity system of second-grade fluid
\begin{equation}\label{1.06}
   \left\{
    \begin{array}{l}
      \omega_t - \alpha^2\Delta\omega_t - \nu\Delta\omega + u \cdot\nabla(\omega - \alpha^2\Delta\omega) =0, \  \ \ \  \text{in} \ \ \Omega \times (0, T), \\
      \ \ \ \ \ \ \ \ \ \ \ \ \ \ \ \ \ \nabla \cdot u = 0, \ \ \ \ \ \ \ \ \ \ \ \ \ \ \ \ \ \ \ \ \ \ \ \  \ \ \ \ \ \ \ \ \ \ \ \ \ \ \ \   \text{in} \ \ \Omega \times (0, T), \\
      \ \ \ \ \ \ \ \ \ \ \ \ \ \ \ \ \ \ \ u = 0, \ \ \ \ \ \ \ \ \ \ \ \ \ \ \ \ \ \ \ \ \ \ \ \ \ \ \  \ \ \ \ \ \ \ \ \ \ \ \ \ \ \ \  \text{on} \ \ \partial\Omega \times (0, T).
    \end{array}
    \right.
\end{equation}
To recover the solution of (\ref{1.01}), (\ref{1.02}) and (\ref{1.04}) from (\ref{1.06}) with initial velocity, we first solve stream function $\psi$ via
\begin{equation}\label{1.07}
   \left\{
    \begin{array}{l}
      \Delta \psi = \omega, ~~~\text{in} ~~\Omega,\\
      \frac{\partial \psi}{\partial n} = 0, ~~~\text{on} ~~ \partial\Omega,
    \end{array}n
    \right.
\end{equation}
 and  obtain the velocity field $u = \nabla^{\bot}\psi = (\partial_{x_2}\psi, -\partial_{x_1}\psi )$  by Boit-Savart law in \cite{Majda2001Vorticity}. With
   aid of the idea of Bona and Wu \cite{Bona2010The}, one shows that the limit of the solution to second-grade fluid model (\ref{1.06}) is exactly the solution to the associated Euler equations as $\nu \rightarrow 0, \alpha \rightarrow 0 $ independently under the radical symmetry situation, we obtain the following main result:

\begin{thm}\label{thm1.01}
Assume that $\omega_0(x)$ is a radially symmetric function in the unit disc $D$
\begin{equation}\label{1.08}
  \omega_0(x) = \omega_0(|x|),
\end{equation}
and   $\omega_0$ to be continuous on $[0, 1]$ such that
$$\int_0^1r\omega_0(r)dr = 0.$$
 Let $u_{0}$ be an initial velocity defined by
  \begin{equation}\label{1.09}
 u_0(|x|) = \frac{1}{r^2}\int_0^r\rho\omega_0(\rho)d\rho
 \left(
  \begin{array}{c}
   -x_2 \\
    x_1 \\
    \end{array}
     \right)
\end{equation}
with $r=|x|.$
Let $u^{(\nu, \alpha)}$ be the solution of the  problem (\ref{1.01}), (\ref{1.02}) and (\ref{1.04})  with initial data satisfies \eqref{1.09}, corresponding to $\nu > 0, \alpha > 0$. Let $u^{E}$ be the solution of Euler system(\ref{1.05}) with the same initial data $u_0$ defined by \eqref{1.09} and the no-penetration  boundary condition
$$u^{E} \cdot n = 0.$$
Let $\omega^{(\nu, \alpha)}$ and $\omega^{E}$ be the vorticities related to the velocity fields $u^{(\nu, \alpha)}$ and $u^{E}$, respectively.
Then for any fixed $T > 0$ and all $t \in [0, T]$, the following differences
$$u^{(\nu, \alpha)}(\cdot, t)- u^{E}(\cdot, t) \ \ \ \text{and} \ \ \ \omega^{(\nu, \alpha)}(\cdot, t)- \omega^{E}(\cdot, t) $$
converge to zero as $\nu \rightarrow 0, \alpha \rightarrow 0$ with $\alpha=o(\nu^{\frac{1}{2}}) $, in $L^2(D)$ and uniformly on compact subset of $D$.
\end{thm}

\begin{rem}
The main theorem (\ref{thm1.01}) and the results in \cite{Filho2015Approximation} tell us that we can verify the convergence from second-grade fluid model to Euler equations as $\alpha \rightarrow 0, \nu \rightarrow 0$, independently. We also find that the limit behavior of second-grade fluid as $\alpha \rightarrow 0, \nu \rightarrow 0$ with $\alpha = o(\nu^{\frac{1}{2}})$ is similar to the limit vanish viscosity from Navier-Stokes equations to Euler under no-slip boundary conditions in radially symmetric case.
\end{rem}

\section{Limit behavior with  radially symmetric}


From the stream function $\psi$ defined by (\ref{1.07}) which is likewise radially symmetric and  rotational invariance of the Laplacian, it is easy to see that  $\psi$ is given as the solution of the ordinary differential equation
\begin{equation}\label{2.01}
  \psi''(r, t) + \frac{1}{r}\psi'(r, t) = \Delta\psi = \omega^{(\nu,\alpha)}(r, t).
\end{equation}
On the other hand, under radial symmetry case, the velocity becomes
\begin{equation}\label{2.02}
  u^{(\nu,\alpha)}(x, t) = (-\partial_{x_2}\psi, \partial_{x_1}\psi ) = (-\frac{x_2}{r}, \frac{x_1}{r})\psi'(r, t),
\end{equation}
where $'$ denotes differentiation with respect to $r$ and $\psi$ satisfies (\ref{2.01}). we can  further reduce it to
$$\frac{d}{dr}(r\psi'(r, t)) = r\omega^{(\nu, \alpha)}(r, t),$$
then it obtains
$$\psi'(r, t)  = \frac{1}{r}\int_0^r\rho\omega^{(\nu, \alpha)}(\rho, t)d\rho.$$
Then, from  the formula (\ref{2.02}) and the above result, one yields
\begin{equation}\label{2.03}
 u^{(\nu,\alpha)}(x, t) = \frac{1}{r^2}\int_0^r\rho\omega^{(\nu, \alpha)}(\rho, t)d\rho
 \left(
  \begin{array}{c}
   -x_2 \\
    x_1 \\
    \end{array}
     \right).
\end{equation}
So, the no-slip boundary conditions $u(x, t) = 0$ imposed when $|x| = 1$ become simply
\begin{equation}\label{2.04}
  \int_0^1r\omega^{(\nu, \alpha)}(r, t)dr = 0.
\end{equation}
Characterized by radial symmetry, it is easy to obtain from the representation (\ref{2.03}) for $u$ that
$$u^{(\nu, \alpha)} \cdot \nabla\omega^{(\nu, \alpha)} = 0,$$
the first  equation of $(\ref{1.06})$ can be reduced to the linear  equation
\begin{equation}\label{2.05}
  \omega_t^{(\nu, \alpha)} - \alpha^2\Delta\omega_t^{(\nu, \alpha)} - \nu\Delta\omega^{(\nu, \alpha)}  =0.
\end{equation}
Similarly,  Euler system (\ref{1.05}) becomes  $\omega^{E}_t= 0$, and thus express the time-independence of the vorticity
\begin{equation}\label{2.06}
  \omega^{E}(x, t) = {\omega}_0(x) = \omega_0(r).
\end{equation}
In consequence, the limit behavior boils down to comparing the solution $\omega^{(\nu, \alpha)}$ of (\ref{2.05}) with initial-boundary value (\ref{1.08}) and (\ref{2.04}). Assume a solution of (\ref{2.05}) in the form
\begin{equation}\label{2.07}
  \omega^{(\nu, \alpha)}(r, t) = Q(r)P(t).
\end{equation}
Substituting (\ref{2.07}) into (\ref{2.05}), which gives
\begin{equation}\label{2.007}
  Q(r)P'(t) - \alpha^2(Q''(r) + \frac{1}{r}Q'(r))P'(t) - \nu(Q''(r) + \frac{1}{r}Q'(r))P(t) = 0.
\end{equation}
There exists a constant $\lambda$ such that
\begin{equation}\label{2.08}
\frac{P'(t)}{P(t)} = \frac{\nu(Q''(r) + \frac{1}{r}Q'(r))}{Q(r) - \alpha^2(Q''(r) + \frac{1}{r}Q'(r))} = -\lambda^2,
\end{equation}
it is equivalent to
\begin{equation}\label{2.09}
  P'(t) + \lambda^2P(t) = 0,
\end{equation}
\begin{equation}\label{2.10}
  Q''(r) + \frac{1}{r}Q'(r) + \frac{\lambda^2}{\nu - \lambda^2\alpha^2}Q(r) = 0.
\end{equation}
Obviously, the solution of (\ref{2.09}) can be given by
\begin{equation}\label{2.11}
  P(t) = P(0)e^{-\lambda^2t}.
\end{equation}
As for equation (\ref{2.10}), it is more subtle. We discuss it in three cases.\\
1) If $\nu = \lambda^2\alpha^2$. \\
     The  equation (\ref{2.007})  can be uniquely solved by  $Q(r) = 0$. By Boit-Savart law $u^{(\nu,\alpha)} \equiv 0$.  In this case, the equation (\ref{2.05})  has a trivial solution. Thus it makes no sense for $\omega(x) = \omega(|x|) \neq 0$.\\
2)  If $\nu < \lambda^2\alpha^2$.\\
   By solving the Bessel equation, it obtains the solution of equation (\ref{2.10}), which is equivalent to equation (\ref{2.007}),
    $$Q(r) = C\sum_{m = 0}^{\infty}\frac{\big(\frac{\lambda}{\sqrt{\lambda^2\alpha^2 - \nu}}r\big)^{2m}}{2^{2m}(m!)^2},$$
    here $C$ is some positive constant. Thus, the equation (\ref{2.05}) has a unique solution
    $$\omega^{(\nu, \alpha)}(r, t) = A\sum_{m = 0}^{\infty}\frac{\big(\frac{\lambda}{\sqrt{\lambda^2\alpha^2 - \nu}}r\big)^{2m}}{2^{2m}(m!)^2}e^{-\lambda^2t} > 0,$$
    here $A = CP(0)$. In this case, the solution $\omega^{(\nu, \alpha)}(r, t)$ is positive for any $(r, t) \in [0, 1) \times [0, T]$. It is incompatible to the boundary condition (\ref{2.04}).\\
3) If $\nu > \lambda^2\alpha^2$.\\
   Solving the Bessel equation (\ref{2.10}), one has
   \begin{equation}\label{2.12}
     Q(r) = CJ_0\big(\frac{\lambda}{\sqrt{\nu - \lambda^2\alpha^2 }}r\big),
   \end{equation}
   here $J_0 = \sum_{m = 0}^{\infty}\frac{(-1)^{m}}{2^{2m}(m!)^2}\big(\frac{\lambda}{\sqrt{\nu - \lambda^2\alpha^2}}r\big)^{2m} $ is the Bessel function of the first kind of order 0 (see Watson\cite{WATSON}). Hence, the solution of equation (\ref{2.05}) in the separated form
   (\ref{2.07}) are
   \begin{equation}\label{2.13}
     \omega^{(\nu, \alpha)}(r, t) = Ae^{-\lambda^2t}J_0\big(\frac{\lambda}{\sqrt{\nu - \lambda^2\alpha^2 }}r\big).
   \end{equation}
 To satisfy the no-slip boundary conditions (\ref{2.04}), one must require
\begin{equation}\label{2.14}
  \int_0^1rJ_0\big(\frac{\lambda}{\sqrt{\nu - \lambda^2\alpha^2 }}r\big)dr = 0.
\end{equation}
New set $y = \frac{\lambda}{\sqrt{\nu - \lambda^2\alpha^2 }}r$ and use the properties of Bessel functions (see  \cite{WATSON})
$$\frac{d}{dy}yJ_1(y) = yJ_0(y).$$
One has
\begin{eqnarray*}
  0 &=& \int_0^1rJ_0\big(\frac{\lambda}{\sqrt{\nu - \lambda^2\alpha^2 }}r\big)dr \\
   &=&  \frac{\nu - \lambda^2\alpha^2}{\lambda^2}\int_0^{\frac{\lambda}{\sqrt{\nu - \lambda^2\alpha^2 }}}\frac{\lambda}{\sqrt{\nu - \lambda^2\alpha^2 }}rJ_0\big(\frac{\lambda}{\sqrt{\nu - \lambda^2\alpha^2 }}r\big)d\frac{\lambda}{\sqrt{\nu - \lambda^2\alpha^2 }}r \\
   &=& \frac{\sqrt{\nu - \lambda^2\alpha^2 }}{\lambda}J_1\big(\frac{\lambda}{\sqrt{\nu - \lambda^2\alpha^2 }}\big).
\end{eqnarray*}
Consequently, the constant $\lambda$ must be satisfied
\begin{equation}\label{2.15}
  J_1\big(\frac{\lambda}{\sqrt{\nu - \lambda^2\alpha^2 }}\big) = 0,
\end{equation}
i.e. $\frac{\lambda}{\sqrt{\nu - \lambda^2\alpha^2 }}$ is a zero of $J_1$. Let $j_1 < j_2 < j_3 < \cdot\cdot\cdot$ denote the sequence of zero
of $J_1$. Thus (\ref{2.15}) is valid only for the discrete set of modes
$$\lambda_k = j_k\sqrt{\nu - \lambda_k^2\alpha^2}, \ \ \text{i.e.} \ \ \lambda_k = \frac{j_k\sqrt{\nu}}{\sqrt{1 + j_k^2\alpha^2}}$$
and so, the solution of (\ref{2.04}) and (\ref{2.05}) in the separated  form (\ref{2.07}) are
$$A_ke^{-\frac{j_k^2\nu}{1 + j_k^2\alpha^2}t}J_0(j_kr),$$
where $A_k$ is a constant. Since the equation (\ref{2.05}) is linear, one expects the complete set of solutions to be given by the series
\begin{equation}\label{2.16}
  \sum_{k = 0}^{\infty}A_ke^{-\frac{j_k^2\nu}{1 + j_k^2\alpha^2}t}J_0(j_kr).
\end{equation}
Now it suffices to determine the coefficients $A_{k}$, thanks to the theory of Dini expansions that generalize the usual Fourier-Bessel expansions \cite{WATSON}. In fact,  for any function
$h$ in the class  $L^2(0, 1; rdr)$. Define
\begin{equation}\label{2.17}
  A_l = \frac{2}{J_0(j_l)^2}\int_0^1h(r)J_0(j_lr)rdr.
\end{equation}
Hence, $h(r)$ can be expressed
\begin{equation}\label{2.18}
  h(r) = 2\int_0^1rh(r)dr + \sum_{k = 1}^{\infty}A_kJ_0(j_kr),
\end{equation}
with the convergence taking place in $L^2(0, 1; rdr)$  in case $h(r)$ is continuous on $[0, 1)$, is uniformly continuous on compact subsets of $[0, 1)$.

Now, applying (\ref{2.18}) to  $\omega^{\nu, \alpha}(r, t)$, one has
$$ \omega^{(\nu, \alpha)}(r, t) = 2\int_0^1r\omega^{(\nu, \alpha)(r, t)}dr +  \sum_{k = 1}^{\infty}A_kJ_0(j_kr) $$
imposing the boundary condition (\ref{2.04}) and initial condition $\omega^{\nu, \alpha}(r, 0) = \omega_0(r), r \in [0, 1)$, and using the orthogonality relations
$$\frac{2}{J_0(j_k)^2}\int_0^1J_0(j_kr)J_0(j_lr)rdr = \delta_{k, l},$$
$\delta_{k,l}$ is the Kronecker delta function.  Using  (\ref{2.17}) on $\omega_0$ which leads to (\ref{2.16}) where
$$A_k = \frac{2}{J_0(j_k)^2}\int_0^1\omega_0(r)J_0(j_kr)rdr.$$

It is now demonstrated that the conditions that $\omega_0(r)$ is continuous in $[0, 1)$ and
$$ \int_0^1r\omega_0(r)dr = 0$$
are sufficient to infer the following limit
\begin{equation}\label{2.19}
  \omega^{(\nu, \alpha)}(r, t) \rightarrow \omega_0(r),  \ \ \ \ \text{as} \ \ \nu \rightarrow 0, \ \alpha \rightarrow 0,
\end{equation}
for $r \in (0, 1)$  and uniformly for $t \in [0, T]$.

Indeed,  in the case $\frac{\alpha^{2}}{\nu}\to 0$, as $\alpha,\nu\to0$, given an arbitrary positive number $\epsilon$ and a compact set $K \subset [0, 1)$, there is an $N$ such that for $n \geq N$,
$$\bigg|\sum_{k = N + 1}^{\infty}A_ke^{-\frac{j_k^2\nu}{1 + j_k^2\alpha^2}t}J_0(j_kr)\bigg|< \bigg|\sum_{k = N + 1}^{\infty}A_kJ_0(j_kr)\bigg| < \frac{\epsilon}{4}$$
for $r \in K$ and $0 \leq t \leq T$. Once $N$ is fixed, we have the fact
$$\frac{j_k^2\nu}{1 + j_k^2\alpha^2}   \rightarrow 0, \ \ \ \ \text{as} \ \ \nu \rightarrow 0, \ \alpha \rightarrow 0,~\mbox{for any}~k=1,2,\cdots N,$$
thus, one yields
$$\bigg|\sum_{k = 1}^N\big(1 - e^{-\frac{j_k^2\nu}{1 + j_k^2\alpha^2}t}\big)A_kJ_0(j_kr)\bigg| < \frac{\epsilon}{2}.$$
Therefor, any $r \in K$ and $t \in [0, T]$, it follows that
\begin{eqnarray*}
  \bigg|\omega^{(\nu, \alpha)}(r, t)- \omega^{E}(r, t)\bigg| &=& \bigg|\sum_{k = 1}^{\infty}A_ke^{-\frac{j_k^2\nu}{1 + j_k^2\alpha^2}t}J_0(j_kr)
   - \sum_{k = 1}^{\infty}A_kJ_0(j_kr)\bigg| \\
   &\leq&  \frac{\epsilon}{2} + \frac{\epsilon}{4} + \frac{\epsilon}{4}=\epsilon.
\end{eqnarray*}
The pointwise convergence of $\omega^{(\nu, \alpha)}$ to $\omega^{E}$ for $x \in D$ immediately implies that the velocity field does. According
to (\ref{2.03}),
$$u^{(\nu, \alpha)}(x, t) - u^{E}(x, t) = \frac{1}{r^2}\int_0^r\rho[\omega^{(\nu, \alpha)}(\rho, t)- \omega^{E}(\rho, t)]d\rho\left(
                                                                                                                                          \begin{array}{c}
                                                                                                                                            -x_2 \\
                                                                                                                                            x_1 \\
                                                                                                                                          \end{array}
                                                                                                                                        \right).
 $$
It follows from the Dominated Convergence Theorem that for any $t \in [0, T]$, $u^{(\nu, \alpha)}(x, t) \rightarrow u^{E}(x, t)$ as $\nu \rightarrow 0, \alpha \rightarrow 0$ for any $x \in D$. The convergence of $u^{(\nu, \alpha)}(x, t)$ to $u^{E}(x, t)$ is uniform on compact subset of $D$
and on the interval $[0, T]$.

It is now proposed that the pointwise convergence implies convergence of both the vorticity and velocity field in $L^2$, which is to say
$$\omega^{(\nu, \alpha)}(x, t) \rightarrow \omega^{E}(x, t) \ \ \ \text{and} \ \ \  u^{(\nu, \alpha)}(x, t) \rightarrow u^{E}(x, t) \ \ \ \text{in} \ \ L^2(D)$$
as $\nu \rightarrow 0, \alpha \rightarrow 0$, uniformly for $t \in [0, T]$ for any fixed $T > 0$. For smooth initial vorticity $\omega_0$, the vorticity $\omega^{E}$ will remain smooth and, in particular, in $L^2$. Therefore
$$\|\omega^{(\nu, \alpha)}(\cdot, t) - \omega^{E}(\cdot, t)\|_{L^2}^2 = \int_{D}|\omega^{(\nu, \alpha)}(x, t) - \omega^{E}(x, t)|^2dx$$
approaches 0 as $\nu \rightarrow 0, \alpha \rightarrow 0$ by the Dominated Convergence Theorem and the pointwise convergence of $\omega^{(\nu, \alpha)}(x, t)$ and $\omega^{E}(x, t)$ for any $x \in D$.

Immediately, one has got  from (\ref{2.03}) that
$$\|u^{(\nu, \alpha)}(\cdot, t) -u^{E}(\cdot, t)\|_{L^2}^2 = 2\pi\int_0^1\big[\int_0^r\rho(\omega^{(\nu, \alpha)}(\rho, t) - \omega^{E}(\rho, t))d\rho\big]^2\frac{dr}{r},$$
also goes to zero. This completes the proof of main theorem.

\section*{Acknowledgment}
 The research of Wenhuo Su  is partially supported by the National Natural Science Foundation of China(No.11801495) and the Science and Technology Project of  Education Department in Jiangxi Province(No.GJJ180833).  Aibin Zang is partially supported by the National Natural Science Foundation of China(No.11771382) and the Science and Technology Project of  Education Department in Jiangxi Province(No.GJJ170891).


\end{document}